\documentclass[12pt,twoside]{article}
\usepackage[english]{babel}
\usepackage[latin1]{inputenc}
\usepackage{amsmath}
\usepackage{amssymb,amsfonts}
\usepackage{graphicx,calc}
\usepackage{times,amssymb,amscd}

\newtheorem{thm}{Theorem}[section]

\newtheorem{lem}[thm]{Lemma}

\newtheorem{defn}[thm]{Definition}

\newtheorem{rem}[thm]{Remark}
\numberwithin{equation}{section}

\newcommand{\bG}{\mathbf{G}}
\newcommand{\bH}{\mathbf{H}}

\newcommand{\bL}{\mathbf{L}}

\newcommand{\bR}{\mathbf{R}}
\newcommand{\bS}{\mathbf{S}}

\newcommand{\bT}{\mathbf{T}}

\newcommand{\bu}{\mathbf{u}}
\newcommand{\bv}{\mathbf{v}}

\newcommand{\cP}{\mathcal{P}}
\newcommand{\cS}{\mathcal{S}}

\newcommand{\cH}{\mathcal{H}}

\newcommand{\EUC}{\mathbb E^3}

\newcommand{\SXR}{\bS^2\!\times\!\bR}
\newcommand{\HXR}{\bH^2\!\times\!\bR}

\newcommand{\SLR}{\widetilde{{\rm SL}_2(\mathbb R)}}
\newcommand{\NIL}{\mathbf{Nil}}
\newcommand{\SOL}{\mathbf{Sol}}

\begin{document}
\pagestyle{myheadings}
\markboth{\centerline{G\'eza Csima and Jen\H o Szirmai}}
{The sum of the interior angles in geodesic and translation...}
\title
{The sum of the interior angles in geodesic and translation triangles of $\SLR$ geometry 
}

\author{G\'eza Csima and Jen\H o Szirmai \\
\normalsize Budapest University of Technology and \\
\normalsize Economics Institute of Mathematics, \\
\normalsize Department of Geometry \\
\normalsize Budapest, P. O. Box: 91, H-1521 \\
\normalsize szirmai@math.bme.hu
\date{\normalsize{\today}}}

\maketitle
\begin{abstract}
We study the interior angle sums of translation and geodesic triangles in the universal cover of $\SLR$ geometry.
We prove that the angle sum $\sum_{i=1}^3(\alpha_i) \ge \pi$ for translation triangles and for geodesic triangles the angle sum can 
be larger, equal or less
than $\pi$.
\end{abstract}
%\tableofcontents

\section{Introduction} \label{section1}
In this paper we are interested in {\it geodesic and translation} triangles in $\SLR$ space that is one of the eight Thurston geometries 
\cite{S,T}. This twisted space can be derived from the 3-dimensional Lie 
group of all $2\times 2$ real matrices with unit determinant.
The space of left invariant Riemannian metrics on the group $\SLR$ is $3$-dimensional \cite{Mi}.
In Section \ref{section2} we describe the projective model of $\SLR$ and we shall use its standard Riemannian metric obtained by 
pull back transform to the infinitesimal  arc-length-square at the origin, 
coinciding also with the global projective metric belonging to the quadratic form (\ref{2.2}). 
We also describe the isometry group of $\SLR$.  In Section \ref{section3} we give an overview not only about geodesic, 
but also about translation curves summarizing previous results, which are related to the projective model of $\SLR$ suggested and introduced 
in \cite{M97}. Our main results will be presented in Section \ref{section4}, namely the possible sum of the interior angles in 
a translation triangle must be greater or equal than $\pi$. However, in geodesic triangles this sum is 
less, greater or equal to $\pi$. 
\section{The projective model for $\SLR$} \label{section2}
\begin{figure}[ht]
\centering
\includegraphics[width=12cm]{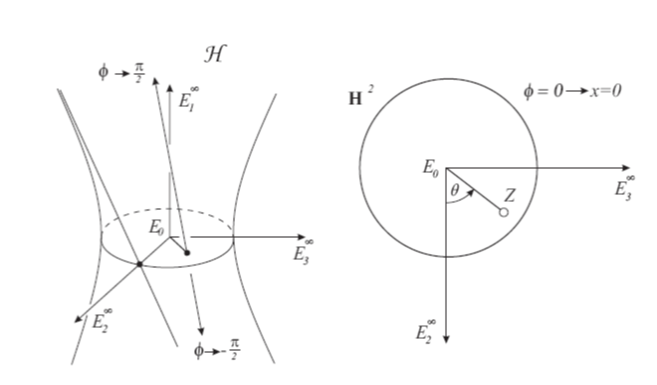}
\caption{The hyperboloid model.}
\label{model}
\end{figure}

Real $ 2\times 2$ matrices $\begin{pmatrix}
         d&b \\
         c&a \\
         \end{pmatrix}$ with the unit determinant $ad-bc=1$
constitute a Lie transformation group by the standard product operation, taken to act on row matrices as point coordinates
\begin{equation}
(z^0,z^1)\begin{pmatrix}
         d&b \\
         c&a \\
         \end{pmatrix}=(z^0d+z^1c, z^0 b+z^1a)=(w^0,w^1)  \label{2.1}
\end{equation}
with
$$
w=\frac{w^1}{w^0}=\frac{\displaystyle b+ (z^1 / z^0) \, a}{\displaystyle d+ ( z^1 / z^0) \, c}=\frac{b+za}{d+zc},
$$
$z=z^{1}/z^{0}$, on the complex projective line ${\mathbb C}^\infty$. This group is a $3$-dimensional manifold, because of its $3$ independent real coordinates and with its usual neighborhood topology \cite{S, T}. In order to model the above structure in the projective sphere $\cP \cS^3$ and in the projective space $\cP^3$ (see \cite{M97}), we introduce the new projective coordinates $(x^0:x^1:x^2:x^3)$ where
$$
a:=x^0+x^3, \quad b:=x^1+x^2, \quad c:=-x^1+x^2, \quad d:=x^0-x^3,
$$
with positive, then the non-zero multiplicative equivalence as a projective freedom in $\cP \cS^3$ and in $\cP^3$, respectively. Then it follows that
\begin{equation}
0>bc-ad=-x^0x^0-x^1x^1+x^2x^2+x^3x^3 \label{2.2}
\end{equation}
describes the internal of the %above
one-sheeted hyperboloid solid $\cH$ in the usual Euclidean coordinate simplex, with the origin $E_0(1:0:0:0)$ and the ideal points of the axes $E_1^\infty(0:1:0:0)$, $E_2^\infty(0:0:1:0)$, $E_3^\infty(0:0:0:1)$. We consider the collineation group ${\bf G}_*$ that acts on the projective sphere $\cP \cS^3$  and preserves the polarity, i.~e. a scalar product of signature $(- - + +)$. This group leaves the one-sheeted hyperboloid solid $\cH$ invariant. We have to choose an appropriate subgroup $\mathbf{G}$ of $\mathbf{G}_*$ as isometry group, then the universal covering group and space $\widetilde{\cH}$ of $\cH$ will be the hyperboloid model of $\SLR$ (see Figure \ref{model} and \cite{M97}).
Consider isometries  given by matrices
\begin{equation}
\begin{gathered}
S(\phi)=
\begin{pmatrix}
\cos{\phi}&\sin{\phi}&0&0 \\
-\sin{\phi}&\cos{\phi}&0&0 \\
0&0&\cos{\phi}&-\sin{\phi} \\
0&0&\sin{\phi}&\cos{\phi}
\end{pmatrix},
\end{gathered} \label{2.3}
\end{equation}
where $\phi \in [0; 2\pi)$. They constitute a one parameter 
group which we denote by $\bS(\phi)$. The elements of $\bS(\phi)$ are the so-called \emph{fibre translations}. 
We obtain a unique fibre line to each $X(x^0:x^1:x^2:x^3) \in \widetilde{\cH}$ as the orbit by right action of 
$\bS(\phi)$ on $X$. The coordinates of points, lying on the fibre line through $X$, can be expressed as the images of $X$ by $\bS(\phi)$:
\begin{equation}
\begin{gathered}
(x^0:x^1:x^2:x^3) \stackrel{\bS(\phi)}{\longrightarrow} {(x^0 \cos{\phi}-x^1 \sin{\phi}: 
x^0 \sin{\phi} + x^1 \cos{\phi}:} \\ {x^2 \cos{\phi} + x^3 \sin{\phi}:-x^2 \sin{\phi}+ x^3 \cos{\phi})}.
\end{gathered} \label{2.4}
\end{equation}
In (\ref{2.3}) and (\ref{2.4}) we can see the $2\pi$-periodicity by $\phi$. 
It admits us to introduce the extension to $\phi \in \mathbb R$, as real parameter, to obtain 
the universal covers $\widetilde{\cH}$ and $\SLR$, respectively, through the projective sphere $\cP\cS^3$. 
The elements of the isometry group of  ${\rm SL}_2(\mathbb R)$ (and so by the above extension the isometries of 
$\SLR$) can be described by the matrix $(a_i^j)$ (see \cite{M97,MSz})
\begin{equation}
 (a_i^j)=
\begin{pmatrix}
a_0^0&a_0^1&a_0^2&a_0^3 \\
\mp a_0^1&\pm a_0^0&\pm a_0^3&\mp a_0^2 \\
a_2^0&a_2^1&a_2^2&a_2^3 \\
\pm a_2^1&\mp a_2^0&\mp a_2^3&\pm a_2^2
\end{pmatrix} , \label{2.5}
\end{equation}
where
$$
\begin{cases} \begin{array}{rcl}
- (a_0^0)^2-(a_0^1)^2+(a_0^2)^2+(a_0^3)^2 & =  & -1, \cr
 - (a_2^0)^2-(a_2^1)^2+(a_2^2)^2+(a_2^3)^2 & =&  1, \cr
- a_0^0a_2^0-a_0^1a_2^1+a_0^2a_2^2+a_0^3a_2^3 & = &  0, \cr
- a_0^0a_2^1+a_0^1a_2^0-a_0^2a_2^3+a_0^3a_2^2 & = & 0,
\end{array}
\end{cases}
$$
and we allow positive proportionality, of course, as projective freedom.

We define the \emph{translation group} $\bG_T$, as a subgroup of the isometry group of ${\rm SL}_2(\mathbb R)$, those isometries acting transitively on the points of ${\cH}$ and by the above extension on the points of $\widetilde{\cH}$. $\bG_T$ maps the origin $E_0(1:0:0:0)$ onto $X(x^0:x^1:x^2:x^3) \in {\cH}$. These isometries and their inverses (up to a positive determinant factor) can be given by
\begin{equation}
\bT=
\begin{pmatrix}
x^0&x^1&x^2&x^3 \\
-x^1&x^0&x^3&-x^2 \\
x^2&x^3&x^0&x^1 \\
x^3&-x^2&-x^1&x^0
\end{pmatrix},\quad \bT^{-1}=
\begin{pmatrix}
x^0&-x^1&-x^2&-x^3 \\
x^1&x^0&-x^3&x^2 \\
-x^2&-x^3&x^0&-x^1 \\
-x^3&x^2&x^1&x^0
\end{pmatrix}. \label{2.6}
\end{equation}

Horizontal intersection of the hyperboloid solid $\cH$  with the plane $E_0 E_2^\infty E_3^\infty$ provides the {\it base plane} of the model $\widetilde{\cH}=\SLR$.

We generally introduce a so-called hyperboloid parametrization by \cite{M97} as follows
\begin{equation}
\begin{cases}
\begin{array}{l}
\displaystyle x^0=\cosh{r} \cos{\phi}, \\
\displaystyle x^1=\cosh{r} \sin{\phi}, \\
\displaystyle x^2=\sinh{r} \cos{(\theta-\phi)}, \\
\displaystyle x^3=\sinh{r} \sin{(\theta-\phi)},
\end{array}
\end{cases}
\label{2.7}
\end{equation}
where $(r,\theta)$ are the polar coordinates of the $\mathbf{H}^2$ base plane, and $\phi$ is the fibre coordinate. We note that
$$
-x^0x^0-x^1x^1+x^2x^2+x^3x^3=-\cosh^2{r}+\sinh^2{r}=-1<0.
$$
The inhomogeneous coordinates, which will play an important role in the later $\EUC$-visualization of the prism tilings in $\SLR$, are given by
\begin{equation} \begin{cases}
\begin{array}{l}
\displaystyle x:=\frac{x^1}{x^0}=\tan{\phi}, \\
\displaystyle y:=\frac{x^2}{x^0}=\tanh{r} \, \frac{\cos{(\theta-\phi)}}{\cos{\phi}}, \\
\displaystyle z:=\frac{x^3}{x^0}=\tanh{r} \, \frac{\sin{(\theta-\phi)}}{\cos{\phi}}.
\end{array}  \end{cases} \label{2.8}
\end{equation}
\section{Geodesic and translation curves} \label{section3}
The infinitesimal arc-length-square of $\SLR$ can be derived by the standard pull back method. By $\bT^{-1}$-action presented by (\ref{2.6}) on differentials $(\mathrm{d}x^0:\mathrm{d}x^1:\mathrm{d}x^2:\mathrm{d}x^3)$,  we obtain the infinitesimal arc-length-square at any point of $\SLR$ in coordinates $(r, \theta, \phi)$:
$$
(\mathrm{d}s)^2=(\mathrm{d}r)^2+\cosh^2{r} \sinh^2{r}(\mathrm{d}\theta)^2+\big[(\mathrm{d}\phi)+\sinh^2{r}(\mathrm{d}\theta)\big]^2.
$$

Hence we get the symmetric metric tensor field $g_{ij}$ on $\SLR$ by components:
\begin{equation}
g_{ij}^*:=
       \begin{pmatrix}
         1&0&0 \\
         0&\sinh^2{r}(\sinh^2{r}+\cosh^2{r})& \sinh^2{r} \\
         0&\sinh^2{r}&1 \\
         \end{pmatrix}.
 \end{equation}
\begin{rem} Similarly to the above computations we obtain the metric tensor for coordinates $(x^1,x^2,x^3)$:
\begin{equation}
g_{ij}:=
\begin{pmatrix}
\frac{1+(x^2)^2+(x^3)^2}{(-1-(x^1)^2+(x^2)^2+(x^3)^2)^2}&\frac{-x^1x^2-2x^3}{(-1-(x^1)^2+(x^2)^2+(x^3)^2)^2}&\frac{-x^1x^3+2x^2)}{(-1-(x^1)^2+(x^2)^2+(x^3)^2)^2} \\
\frac{-x^1x^2-2x^3}{(-1-(x^1)^2+(x^2)^2+(x^3)^2)^2}&\frac{1+(x^1)^2+(x^3)^2}{(-1-(x^1)^2+(x^2)^2+(x^3)^2)^2}& \frac{x^2x^3}{(-1-(x^1)^2+(x^2)^2+(x^3)^2)^2} \\
\frac{-x^1x^3+2x^2}{(-1-(x^1)^2+(x^2)^2+(x^3)^2)^2}&\frac{x^2x^3}{(-1-(x^1)^2+(x^2)^2+(x^3)^2)^2}&\frac{1+(x^1)^2+(x^2)^2}{(-1-(x^1)^2+(x^2)^2+(x^3)^2)^2}\\
         \end{pmatrix}. \label{3.1}
 \end{equation}
\end{rem}
\subsection{Geodesic curves}
The geodesic curves of $\SLR$ are generally defined as having locally minimal arc length between any two of their (close enough) points.

By (\ref{3.1}) the second order differential equation system of the $\SLR$ geodesic curve is the following:
\begin{equation}
\begin{cases}
\begin{array}{l}
\displaystyle
\ddot{r}=\sinh(2r)~\! \dot{\theta}~\! \dot{\phi}+\frac12 \big[ \sinh(4r)-\sinh(2r) \big] \dot{\theta} ~\! \dot{\theta},\\
\displaystyle  \ddot{\theta}=\frac{2\dot{r}}{\sinh{(2r)}} \big[ (3 \cosh{(2r)}-1) \dot{\theta}+2\dot{\phi} \big], \\
\displaystyle  \ddot{\phi} = 2\dot{r}\tanh{(r)} \, [2\sinh^2{(r)}~\! \dot{\theta}+ \dot{\phi}].
\end{array}
\end{cases} \label{2.10}
\end{equation}

We can assume, by the homogeneity, that the starting point of a geodesic curve is the origin $(1:0:0:0)$. Moreover,
$$
\begin{cases}
\begin{array}{l}
r(0) = 0, \\
\theta (0)  = 0,  \\
\phi (0) = 0, \\
\end{array}
\end{cases}
\qquad \mbox{and} \qquad
\begin{cases}
\begin{array}{l}
\dot{r} (0) = \cos(\alpha), \\
\dot{\theta}(0)= - \sin(\alpha), \\
\dot{\phi}(0)=\sin(\alpha)
\end{array}
\end{cases}
$$
are the initial values in Table 1 for the solution of (\ref{2.10}), and so the unit velocity will be achieved (see details in \cite{DESS}). The solutions are parametrized by the arc-length $s$ and the angle $\alpha$ from the initial condition: $(r(s,\alpha), \theta (s, \alpha), \phi (s,\alpha))$.

\begin{table}[ht]
\caption{Geodesic curves.} \label{table1}
\vspace{3mm}
\centerline{
$
\begin{array}{|c|l|} \hline
\textrm{direction} & \textrm{parametrization of a geodesic curve}  \cr \hline
\begin{gathered}
0 \le \alpha < \frac{\pi}{4} \\ (\bH^2-{\rm like})
\end{gathered}
&
\begin{array}{l}
r(s,\alpha) =  {\mathrm{arsinh}} \Big( \frac{\cos{\alpha}}{\sqrt{\cos{2\alpha}}}\sinh(s\sqrt{\cos{2\alpha}}) \Big) \cr
\theta(s,\alpha)  =  -{\mathrm{arctan}} \Big( \frac{\sin{\alpha}}{\sqrt{\cos{2\alpha}}}\tanh(s\sqrt{\cos{2\alpha}}) \Big) \cr
\phi(s,\alpha)  =  2 s\sin{\alpha} + \theta(s,\alpha)
\end{array} \cr \hline
\begin{gathered} \alpha=\frac{\pi}{4} \\ ({\rm light-like}) \end{gathered}
&
\begin{array}{l}
r(s,\alpha)={\mathrm{arsinh}} \Big( \frac{\sqrt{2}}{2} s \Big) \cr
\theta(s,\alpha)=-{\mathrm{arctan}} \Big( \frac{\sqrt{2}}{2} s \Big) \cr
\phi(s,\alpha)=\sqrt{2} s +\theta(s,\alpha)
\end{array} \cr \hline
\begin{gathered} \frac{\pi}{4}  < \alpha \le \frac{\pi}{2} \\ ({\rm fibre-like}) \end{gathered}
&
\begin{array}{l}  r(s,\alpha)={\mathrm{arsinh}} \Big( \frac{\cos{\alpha}}{\sqrt{-\cos{2\alpha}}}\sin(s\sqrt{-\cos{2\alpha}}) \Big) \cr
\theta(s,\alpha)=-{\mathrm{arctan}} \Big( \frac{\sin{\alpha}}{\sqrt{-\cos{2\alpha}}}\tan(s\sqrt{-\cos{2\alpha}}) \Big) \cr
\phi(s,\alpha)=2s\sin{\alpha} + \theta(s,\alpha)
\end{array}  \cr \hline
\end{array}
$
}
\end{table}

The parametrization of a geodesic curve in the hyperboloid model with the geographical sphere coordinates $(\lambda, \alpha)$, as longitude and altitude, $(-\pi < \lambda \le \pi, ~ -\frac{\pi}{2}\le \alpha \le \frac{\pi}{2})$, and the arc-length parameter $s \geq 0$, has been determined in \cite{DESS}. The Euclidean coordinates $X(s,\lambda,\alpha)$, $Y(s,\lambda,\alpha)$, $Z(s,\lambda,\alpha)$ of the geodesic curves can be determined by substituting the results of Table 1 (see also \cite{DESS}) into formula (\ref{2.8}) as follows
\begin{equation}
\begin{cases}
\begin{array}{l}
\displaystyle X(s,\lambda,\alpha)=\tan{(\phi(s,\alpha))}, \\
\displaystyle Y(s,\lambda,\alpha)=\frac{\tanh{(r(s,\alpha))}}{\cos{(\phi(s,\alpha))}} \cos \big[ \theta(s,\alpha)-\phi(s,\alpha)+\lambda \big],\\
\displaystyle Z(s,\lambda,\alpha)=\frac{\tanh{(r(s,\alpha))}}{\cos{(\phi(s,\alpha))}} \sin \big[ \theta(s,\alpha)-\phi(s,\alpha)+\lambda \big].
\end{array}
\end{cases}
\label{2.11}
\end{equation}
\begin{defn}
The \emph{geodesic distance} $d(P,Q)$ between points $P,Q  \in \SLR$ is defined as the arc length of the geodesic curve from $P$ to $Q$.
\end{defn}
\subsection{Translation curves}
We recall some basic facts about translation curves in $\SLR$  following 
\cite{MSz11,MSz, MSzV14}. For any point $X(x^{0}: x^{1}: x^{2}: x^{3}) \in {\mathcal H}$ 
(and later also for points in $\widetilde{\mathcal H}$) the \emph{translation map} from  the origin 
$E_{0} (1:0:0:0)$ to $X$ is defined by the \emph{translation matrix} $\bT$ and its inverse presented in (\ref{2.6}).

Let us consider for a given vector  $(q:u:v:w)$ a curve 
$\mathcal C(t) = (x^{0}(t):x^{1}(t):$ $x^{2}(t): x^{3}(t))$, $t \geq 0$, in ${\mathcal H}$ starting at the origin: $\mathcal C(0) = E_{0}(1:0:0:0)$ and such that
$$
\dot{\mathcal C} (0) = (\dot{x}^{0}(0): \dot{x}^{1}(0): \dot{x}^{2}(0): \dot{x}^{3}(0)) = (q:u:v:w),
$$
where $\dot{\mathcal C}(t)  = (\dot{x}^{0}(t): \dot{x}^{1}(t): \dot{x}^{2}(t): \dot{x}^{3}(t))$
is the tangent vector at any point of the curve. For $t \geq 0$ there exists a matrix
\begin{equation*}
\bT(t) = \left( \begin{array}{cccc}
x^{0}(t) & x^{1}(t) & x^{2}(t) & x^{3}(t) \cr
-x^{1}(t) & x^{0}(t) & x^{3}(t) & -x^{2}(t) \cr
x^{2}(t) & x^{3}(t) & x^{0}(t) & x^{1}(t) \cr
x^{3}(t) & -x^{2}(t) & -x^{1}(t) & x^{0}(t) \cr
\end{array} \right) \label{2.13}
\end{equation*}
which defines the translation from $\mathcal C(0)$ to $\mathcal C (t)$:
\begin{equation*}
\mathcal C(0) \cdot \bT(t) = \mathcal C(t),  \quad t \geqslant 0. \label{2.14}
\end{equation*}
The $t$-parametrized family $\bT(t)$ of translations is used in the following definition.

\begin{defn}
The curve $\mathcal C(t)$, $t \geqslant 0$, is said to be a \emph{translation curve} if
\begin{equation*}
\dot{\mathcal C }(0) \cdot \bT(t) = \dot{\mathcal C}(t),  \quad t \geqslant 0. \label{2.15}
\end{equation*}
\end{defn}

The solution, depending on $(q:u:v:w)$ had been determined in \cite{MSz11}, where it splits into three cases. 

It was observed above that for any $X(x^{0}:x^{1}:x^{2}:x^{3}) \in \widetilde{\mathcal H}$ there 
is a suitable transformation $\bT^{-1}$, given by (\ref{2.6}), which sent $X$ to the origin $E_{0}$ along a translation curve.
\begin{defn} \label{def-distance}
A  \emph{translation distance} $\rho (E_{0},X)$ between the origin $E_{0}(1:0:0:0)$ and the point  $X(1:x:y:z)$ is the length of a 
translation curve connecting them.
\end{defn}
For a given translation curve $\mathcal C = \mathcal C (t)$ the initial unit tangent vector $(u,v,w)$ (in Euclidean coordinates) 
at $E_{0}$ can be presented as
\begin{equation}
u = \sin \alpha, \quad v = \cos \alpha \cos \lambda, \quad w = \cos \alpha \sin \lambda, \label{2.25}
\end{equation}
for some $-\frac{\pi}{2} \leqslant \alpha \leqslant \frac{\pi}{2}$ and $ -\pi < \lambda \leqslant \pi$.  
In $\widetilde{\mathcal H}$ this vector is of length square $-u^{2} + v^{2} + w^{2}  = \cos 2 \alpha$. We 
always can assume that $\mathcal C$ is parametrized by the translation arc-length parameter $t = s \geqslant 0$. 
Then coordinates of a point $X(x, y, z)$ of $\mathcal C$, such that the translation distance between $E_{0}$ and $X$ equals 
$s$, depend on $(\lambda,\alpha,s)$ as geographic coordinates according to the above considered three cases as follows.
\begin{table}[ht]
\caption{Translation  curves.} \label{table2}
\vspace{3mm}
\centerline{
$
\begin{array}{|c|l|} \hline
\textrm{direction} & \textrm{parametrization of a translation curve}  \cr \hline  
\begin{gathered}
0 \le \alpha < \frac{\pi}{4} \\ (\bH^2-{\rm like})
\end{gathered}
&
\begin{array}{l}
\begin{gathered} \noalign{\vskip2pt}\begin{pmatrix} x(s,\alpha,\lambda)\\y(s,\alpha,\lambda)\\z(s,\alpha,\lambda)\end{pmatrix}=  \frac{\tanh (s \sqrt{\cos 2 \alpha})}{\sqrt{\cos 2\alpha}} \begin{pmatrix} \sin \alpha\\
 \cos \alpha \cos \lambda\\
 \cos \alpha \sin \lambda)\end{pmatrix} \end{gathered} \cr \noalign{\vskip2pt}\end{array} \cr
 \hline
\begin{gathered} \alpha=\frac{\pi}{4} \\ ({\rm light-like}) \end{gathered}
&
\begin{array}{l}
\begin{gathered} \noalign{\vskip2pt}\begin{pmatrix} x(s,\alpha,\lambda)\\y(s,\alpha,\lambda)\\z(s,\alpha,\lambda)\end{pmatrix}=  \frac{\sqrt{2} s}{2} 
\begin{pmatrix} 1\\
 \cos \lambda\\
 \sin \lambda)\end{pmatrix} \end{gathered} \cr
\noalign{\vskip2pt}\end{array} \cr \hline
\begin{gathered} \frac{\pi}{4}  < \alpha \le \frac{\pi}{2} \\ ({\rm fibre-like}) \end{gathered}
&
\begin{array}{l}  \begin{gathered} \noalign{\vskip2pt}\begin{pmatrix} x(s,\alpha,\lambda)\\y(s,\alpha,\lambda)\\z(s,\alpha,\lambda)\end{pmatrix}=  \frac{\tan (s \sqrt{-\cos 2 \alpha})}{\sqrt{-\cos 2\alpha}} \begin{pmatrix} \sin \alpha\\
 \cos \alpha \cos \lambda\\
 \cos \alpha \sin \lambda)\end{pmatrix} \end{gathered}   \cr
\noalign{\vskip2pt}
\end{array}  \cr \hline
\end{array}
$
}
\end{table}
\section{Geodesic and translation triangles} \label{section4}
\subsection{Geodesic triangles}
We consider $3$ points $A_1$, $A_2$, $A_3$ in the projective model of $\SLR$ space (see Section 2). 
The {\it geodesic segments} $a_k$ between the points $A_i$ and $A_j$ 
$(i<j,~i,j,k \in \{1,2,3\}, k \ne i,j$ are called sides of the {\it geodesic triangle} $T_g$ with vertices $A_1$, $A_2$, $A_3$.  

In Riemannian geometries the metric tensor (\ref{3.1}) is used to define the angle $\theta$ between two geodesic curves. 
If their tangent vectors in their common point are $\bu$ and $\bv$ and $g_{ij}$ are the components of the metric tensor then
\begin{equation}
\cos(\theta)=\frac{u^i g_{ij} v^j}{\sqrt{u^i g_{ij} u^j~ v^i g_{ij} v^j}} 
\end{equation}
It is clear by the above definition of the angles and by the metric tensor (3.2), that 
the angles are the same as the Euclidean ones at the origin of 
the projective model of $\SLR$ geometry. 

We note here that the angle of two intersecting geodesic curves depend on the orientation of the tangent vectors. We will consider
the {\it interior angles} of the triangles that are denoted at the vertex $A_i$ by $\omega_i$ $(i\in\{1,2,3\})$.

\subsubsection{Fibre-like right angled triangles}
\begin{figure}[ht]
\centering
\includegraphics[width=12cm]{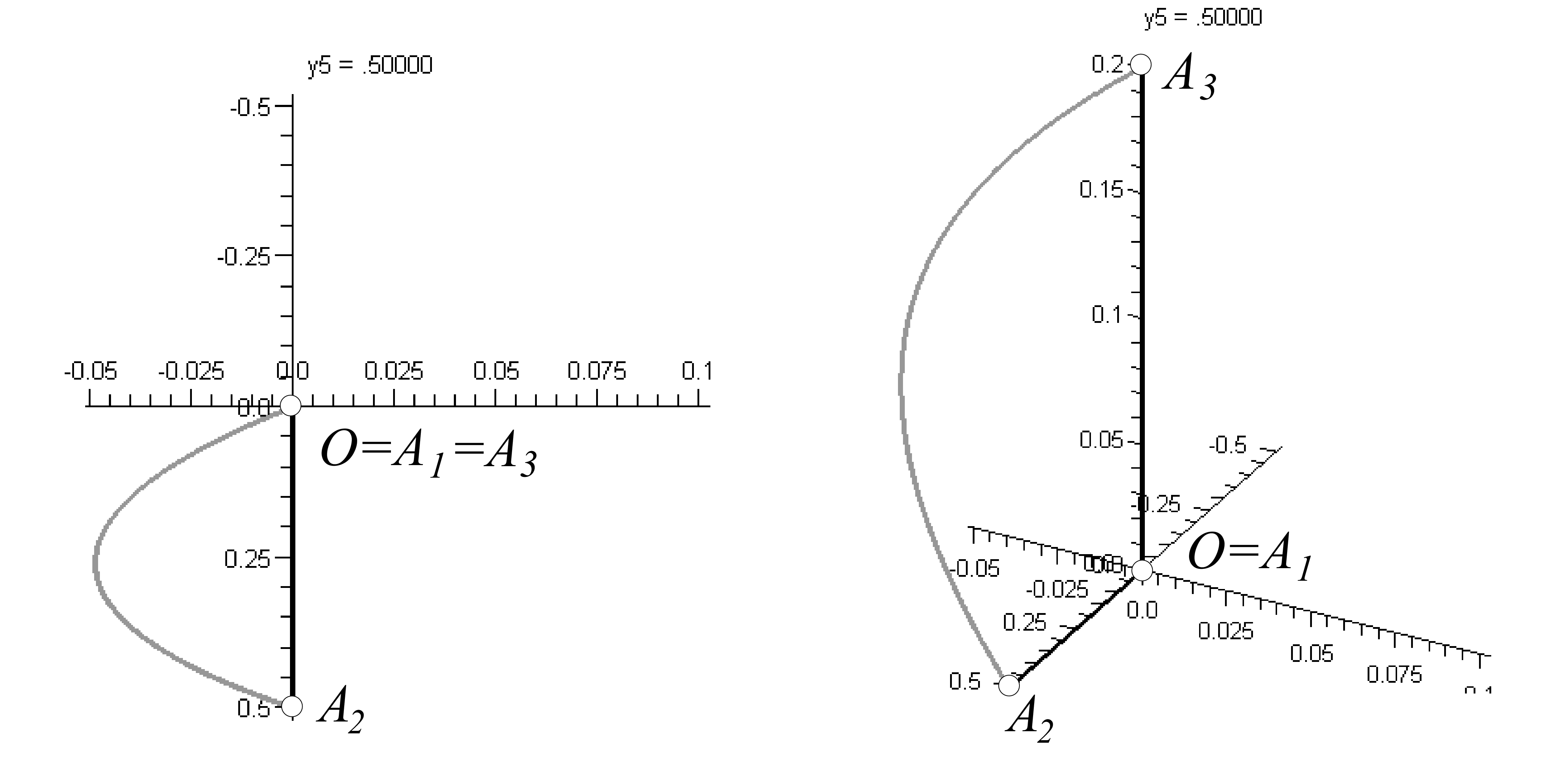}
\caption{Fibre-like geodesic triangle}
\label{fibr}
\end{figure}

A geodesic triangle is called fibre-like if one of its edges lies on a fibre line. In this section we study the right angled fibre-like
triangles. We can assume without loss of generality
that the vertices $A_1$, $A_2$, $A_3$ of a fibre-like right angled triangle (see Figure \ref{fibr}) $T_g$ have the following coordinates:
\begin{equation}
A_1=(1:0:0:0),~A_2=(1:0:y^2:0),~A_3=(1:x^3:0:0) 
\end{equation}
The geodesic segment $A_1A_2$ lies on the $y$ axis, the geodesic segment $A_1A_3$ lies on the $x$ axis (see Table 1) and its angle is
$\omega_1=\frac{\pi}{2}$ in the $\SLR$ space (this angle is in Euclidean sense also $\frac{\pi}{2}$ since $A_1=E_0$).  

\begin{figure}[ht]
\centering
\includegraphics[width=12cm]{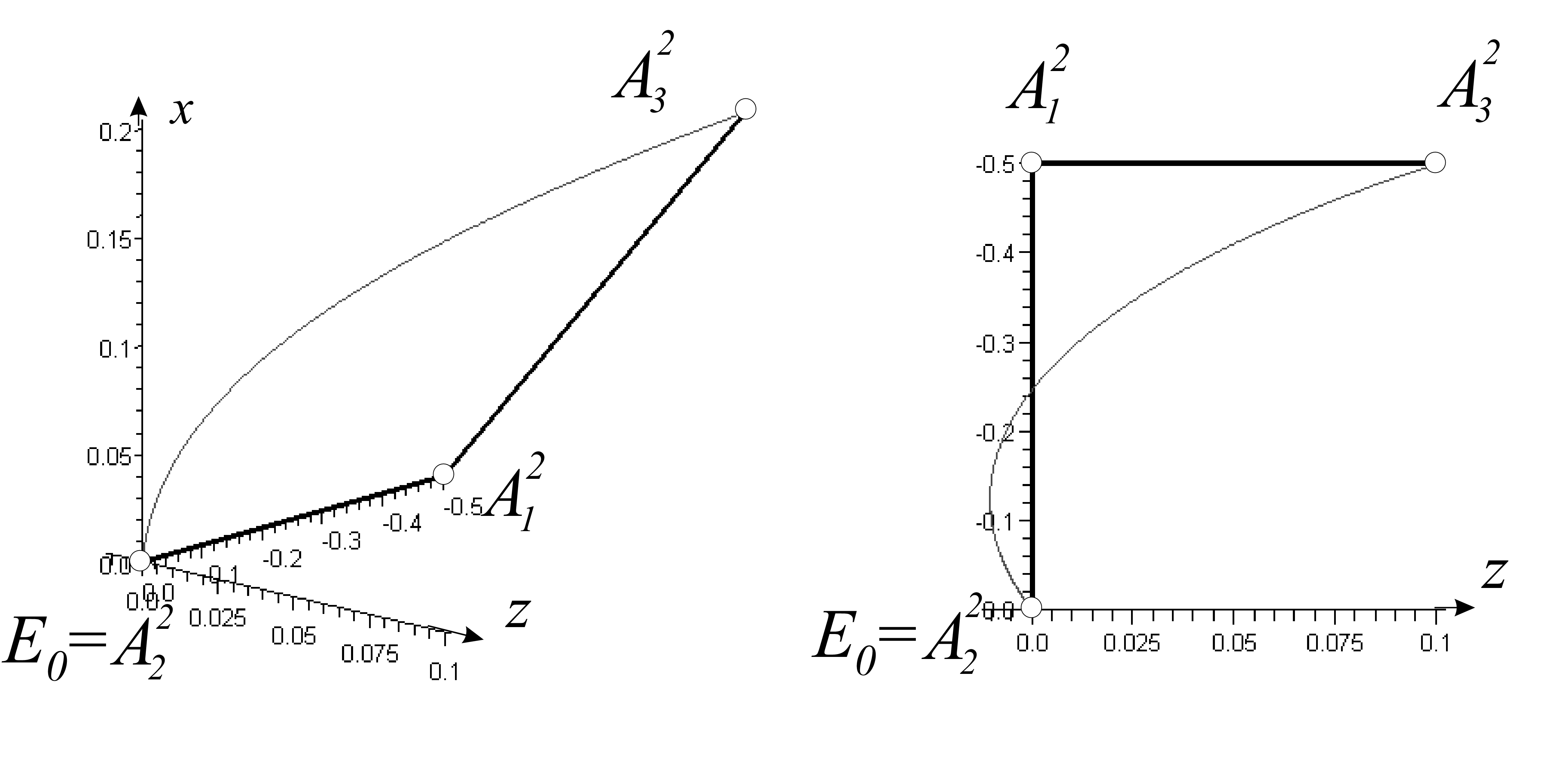}
\caption{Translated image $A_1^2A_2^2A_3^3$ of the fibre-like geodesic triangle $A_1A_2A_3$  by translation $\bT_{A_2}$}
\label{transant}
\end{figure}

In order to determine the further interior angles of fibre-like geodesic triangle $A_1A_2A_3$ 
we define \emph{translations} $\bT_{A_i}$, $(i\in \{2,3\})$ as elements of the isometry group of ${\rm SL}_2(\mathbb R)$, (see \ref{2.6}) that 
maps the origin $E_0$ onto $A_i$. 
E.g. the isometry $\bT_{A_2}$ and its inverse (up to a positive determinant factor) can be given by:

\begin{equation}
\bT_{A_2}=
\begin{pmatrix}
1 & 0 & y^2 & 0 \\
0 & 1 & 0 & -y^2 \\
y^2 & 0 & 1 & 0 \\
0 & -y^2 & 0 & 1
\end{pmatrix} , ~ ~ ~
\bT_{A_2}^{-1}=
\begin{pmatrix}
1 & 0 & -y^2 & 0 \\
0 & 1 & 0 & y^2 \\
-y^2 & 0 & 1 & 0 \\
0 & y^2 & 0 & 1
\end{pmatrix} , 
\end{equation}
and the images $\bT_{A_2}(A_i)$ of the vertices $A_i$ $(i \in \{1,2,3\}$ are the following (see also Figure \ref{transant}):   
\begin{equation}
\begin{gathered}
\bT^{-1}_{A_2}(A_1)=A_1^2=(1:0:-y^2:0),~\bT^{-1}_{A_2}=A_2^2(A_2)=E_0=(1:0:0:0), \\ \bT^{-1}_{A_2}(A_3)=A_3^2=(1:x^3:-y^2:x^3y^2), 
\end{gathered}
\label{4.4}
\end{equation}
Similarly to the above cases we obtain:
\begin{equation}
\begin{gathered}
\bT^{-1}_{A_3}(A_1)=A_1^3=(1:-x^3:0:0),~\bT^{-1}_{A_3}(A_2)=A_2^3=(1:-x^3:y^2:-x^3y^2), \\ \bT^{-1}_{A_3}(A_3)=A_3^3=E_0=(1:0:0:0), 
\end{gathered}
\label{4.5}
\end{equation}
\begin{figure}[ht]
\centering
\includegraphics[width=12cm]{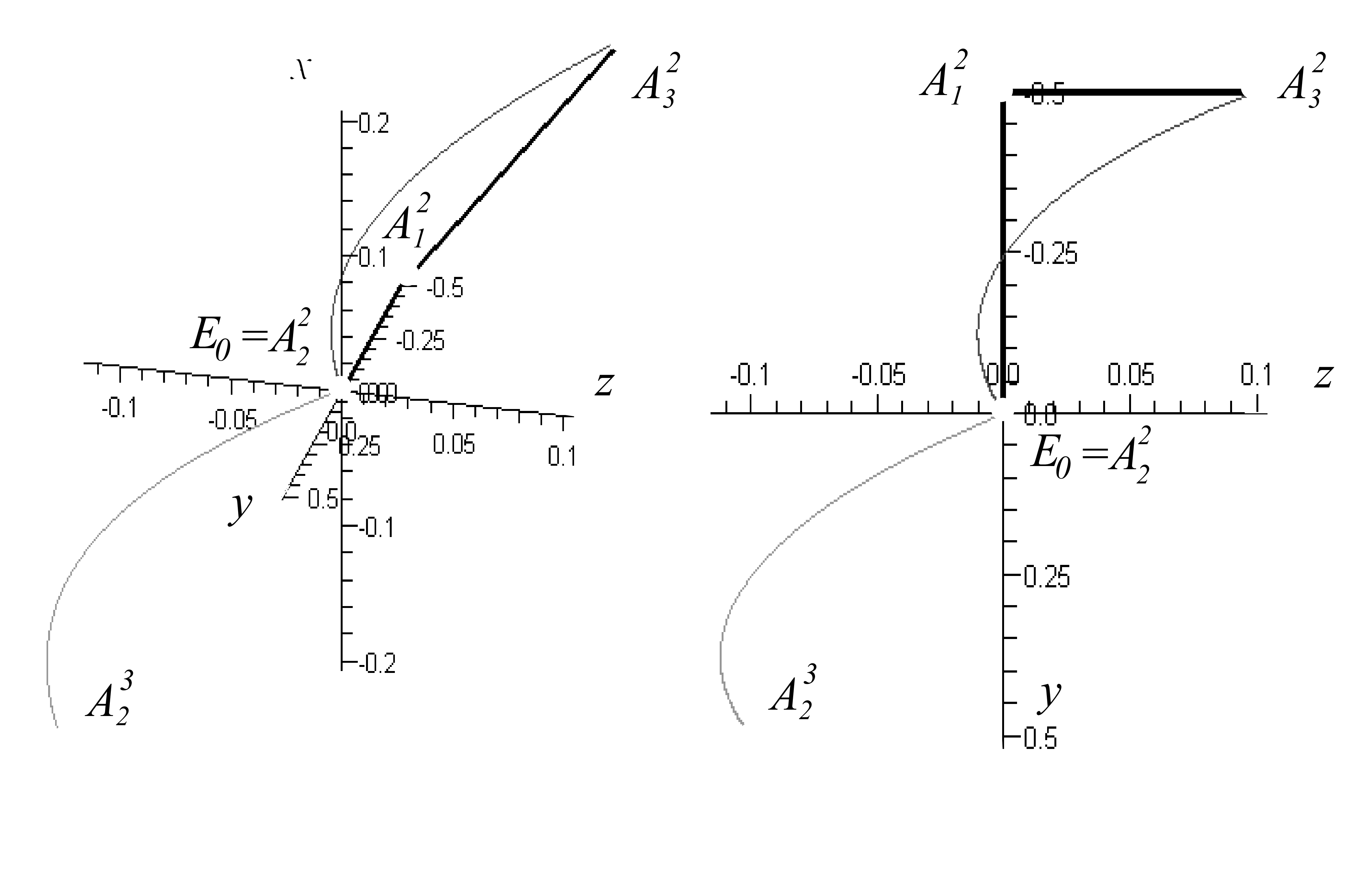}
\caption{The translated points $A_2^3$ and $A_3^2$ are antipodal related to the origin $E_0$}
\label{transfibr}
\end{figure}
Our aim is to determine angle sum $\sum_{i=1}^3(\omega_i)$ of the interior angles of the above right angled 
fibre-like geodesic triangle $A_1A_2A_3$.
We have seen that $\omega_1=\frac{\pi}{2}$ and the angle of geodesic curves with common point at the origin $E_0$ is the same as the 
Euclidean one therefore it can be determined by usual Euclidean sense. Moreover, the translations $\bT_{A_2}$ and $\bT_{A_3}$ are isometries 
in $\SLR$ geometry thus
$\omega_2$ is equal to the angle $(g(A_2^2, A_1^2)g(A_2^2, A_3^2)) \angle$ (see Figure \ref{transfibr}) 
where $g(A_2^2, A_1^2)$, $g(A_2^2, A_3^2)$ are the oriented geodesic curves and $\omega_3$ is equal to the angle 
$(g(A_3^3 A_1^3)g(A_3^3 A_2^3)) \angle$ $(E_0=A_2^2=A_3^3$ see \ref{4.4}, \ref{4.5}).    

The parametrization of a geodesic curve in the model is given by the geographical sphere coordinates 
$(\lambda, \alpha)$, as longitude and altitude, $(-\pi < \lambda \le \pi, ~ -\frac{\pi}{2}\le \alpha \le \frac{\pi}{2})$, 
and the arc-length parameter $s \geq 0$ (see Table 1 and (3.4)). 

We denote the oriented unit tangent vector of the oriented geodesic curves $g(E_0, A_i^j)$ with $\mathbf{t}_i^j$ where
$(i,j)=(2,3),(3,2),(1,2),(1,3)$.
The Euclidean coordinates of $\mathbf{t}_i^j$ are: 
\begin{equation}
\mathbf{t}_i^j=(\sin(\alpha_i^j), \cos(\alpha_i^j) \cos(\lambda_i^j), \cos(\alpha_i^j) \sin(\lambda_i^j)). 
\end{equation}
\begin{lem}
The sum of the interior angles of a fibre-like right angled geodesic triangle is greater or equal to $\pi$. 
\end{lem}
\textbf{Proof:} It is clear, that $\mathbf{t}_1^2=(0,-1,0)$ and $\mathbf{t}_1^3=(-1,0,0)$. Moreover, 
the points $A_2^3$ and $A_3^2$ are antipodal related to the origin $E_0$ therefore the equation $|\alpha_2^3|=|\alpha_3^2|$ holds
(i.e. the angle between the vector $\mathbf{t}_2^3$ and $[y,z]$ plane are equal to the angle between the vector $\mathbf{t}_3^2$ and the $[y,z]$
plane).
That means, that $\omega_3=\frac{\pi}{2}-|\alpha_2^3|=\frac{\pi}{2}-|\alpha_3^2|$.

The vector $\mathbf{t}_1^2$ lies in the $[y,z]$ plane therefore the angle $\omega_2$ is greater or equal than $\alpha_2^3=\alpha_3^2$.
Finally we obtain, that 
$$
\sum_{i=1}^3(\omega_i)=\frac{\pi}{2}+\frac{\pi}{2}-|\alpha_2^3|+\omega_2 \ge \pi. \ \ \square
$$
In the following table we summarize some numerical data of geodesic triangles for given parameters:    
\medbreak
\centerline{\vbox{
\halign{\strut\vrule\quad \hfil $#$ \hfil\quad\vrule
&\quad \hfil $#$ \hfil\quad\vrule &\quad \hfil $#$ \hfil\quad\vrule &\quad \hfil $#$ \hfil\quad\vrule &\quad \hfil $#$ \hfil\quad\vrule &\quad \hfil $#$ \hfil\quad\vrule
\cr
\noalign{\hrule}
\multispan6{\strut\vrule\hfill\bf Table 3, ~ ~  $x_3=1/5$ \hfill\vrule}%
\cr
\noalign{\hrule}
\noalign{\vskip2pt}
\noalign{\hrule}
y^2 & |\alpha_2^3|=|\alpha_3^2| & d(A_2A_3) & \omega_2 & \omega_3  & \sum_{i=1}^3(\omega_i)  \cr
\noalign{\hrule}
\rightarrow 0 & \rightarrow \pi/2  & \begin{gathered} \rightarrow \\ \mathrm{arctan}(1/5) \approx \\
\approx 0.1974 \end{gathered} & \rightarrow \pi/2 & \rightarrow 0 & \rightarrow \pi \cr
\noalign{\hrule}
1/1000 & 1.5657 & 0.1974 & 1.5658 & 0.0051 & 3.1417 \cr
\noalign{\hrule}
1/3 & 0.4993 & 0.3970 & 0.3560 & 1.0715 & 3.1806 \cr
\noalign{\hrule}
1/2 & 0.3170 & 0.5809 & 0.3560 & 1.2538 & 3.1806 \cr
\noalign{\hrule}
3/4 & 0.1630 & 0.9891 & 0.2043 & 1.4078 & 3.1829 \cr
\noalign{\hrule}
999/1000 & 0.0299 & 3.8032 & 0.0422 & 1.5409 & 3.1540\cr
\noalign{\hrule}
\rightarrow 1 & \rightarrow 0 & \rightarrow \infty & \rightarrow 0 & \rightarrow \pi/2& \rightarrow \pi \cr
\noalign{\hrule}
}}}
\medbreak
\subsubsection{Hyperbolic-like right angled geodesic triangles}
\begin{figure}[ht]
\centering
\includegraphics[width=8cm]{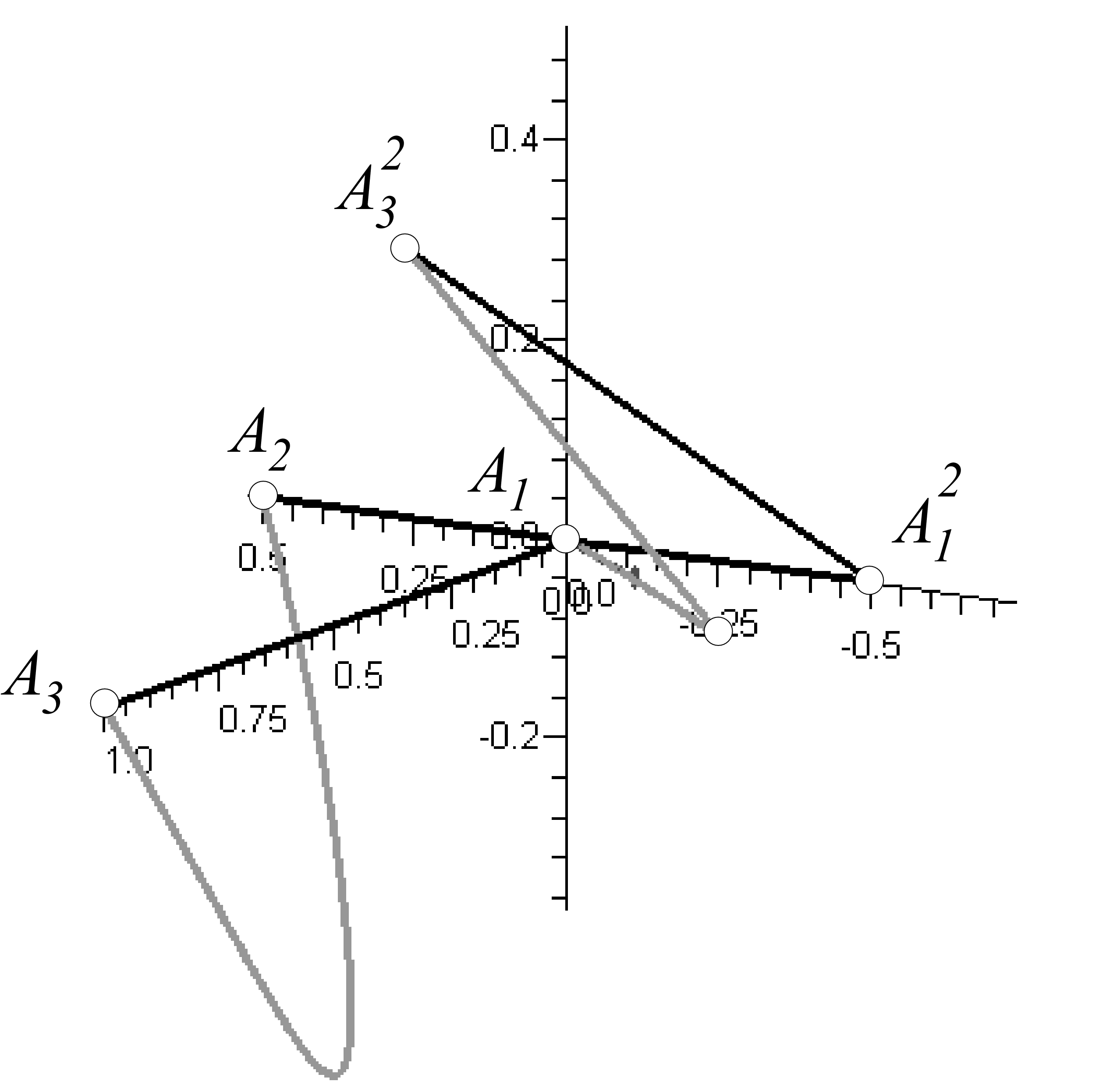}
\caption{Hyperbolic-like geodesic triangle $A_1A_2A_3$ and its translated copy $A_1^2A_2^2A_3^2$.}
\label{hyph}
\end{figure}

A geodesic triangle is hyperbolic-like if its vertices lie in the base plane (i.e. $[y,z]$ coordinate plane) of the model. 
In this section we analyze the interior angle sum of the right angled hyperbolic-like
triangles. We can assume without loss of generality
that the vertices $A_1$, $A_2$, $A_3$ of a hyperbolic-like right angled triangle (see Figure \ref{hyph}) $T_g$ have the following coordinates:
\begin{equation}
A_1=(1:0:0:0),~A_2=(1:0:y^2:0),~A_3=(1:0:0:z^3) 
\end{equation}
The geodesic segment $A_1A_2$ lies on the $y$ axis, the geodesic segment $A_1A_3$ lies on the $z$ axis and its angle is
$\omega_1=\frac{\pi}{2}$ in the $\SLR$ space (this angle is in Euclidean sense also $\frac{\pi}{2}$ since $A_1=E_0$).  

 In order to determine the further interior angles of fibre-like geodesic triangle $A_1A_2A_3$ 
similarly to the fibre-like case we define a \emph{translations} $\bT_{A_i}$, (see \ref{2.6}) that 
maps the origin $E_0$ onto $A_i$. 
E.g. the isometry $\bT_{A_3}$ and their inverses (up to a positive determinant factor) can be given by:
\begin{equation}
\bT_{A_3}=
\begin{pmatrix}
1 & 0 & 0 & z^3 \\
0 & 1 & z^3 & 0 \\
0 & z^3 & 1 & 0 \\
z^3 & 0 & 0 & 1
\end{pmatrix} , ~ ~ ~
\bT_{A_3}^{-1}=
\begin{pmatrix}
1 & 0 & 0 & -z^3 \\
0 & 1 & -z^3 & 0 \\
0 & -z^3 & 1 & 0 \\
-z^3 & 0 & 0 & 1
\end{pmatrix}. 
\end{equation}
We get similarly to the above cases that the images $\bT^{-1}_{A_j}(A_i)$ of the vertices $A_i$ $(i \in \{1,2,3\}, ~ j \in \{2,3\} )$  
are the following (see also Figure \ref{hyph}):   
\begin{equation}
\begin{gathered}
\bT^{-1}_{A_2}(A_1)=A_1^2=(1:0:-y^2:0),~\bT^{-1}_{A_2}=A_2^2(A_2)=E_0=(1:0:0:0), \\ \bT^{-1}_{A_2}(A_3)=A_3^2=(1:y^2z^3:-y^2:z^3), 
\end{gathered}
\end{equation}
\begin{equation}
\begin{gathered}
\bT^{-1}_{A_3}(A_1)=A_1^3=(1:-z^3:0:0),~\bT^{-1}_{A_3}(A_2)=A_2^3=(1:-y^2z^3:y^2:-z^3), \\ \bT^{-1}_{A_3}(A_3)=A_3^3=E_0=(1:0:0:0), 
\end{gathered}
\end{equation}
We study similarly to the above fibre-like case the sum $\sum_{i=1}^3(\omega_i)$ of the interior angles of the above right angled 
hyperbolic-like geodesic triangle $A_1A_2A_3$.

It is clear, that the angle of geodesic curves with common point at the origin $E_0$ is the same as the 
Euclidean one therefore it can be determined by usual Euclidean sense. The translations $\bT_{A_2}$ and $\bT_{A_3}$ preserve 
the measure of angles $\omega_i$ $(i \in \{2,3\}$ therefore $\omega_2=(g(A_2^2, A_1^2)g(A_2^2, A_3^2)) \angle$ (see Figure \ref{hyph}) 
and $\omega_3=(g(A_3^3 A_1^3)g(A_3^3 A_2^3)) \angle$.    

Similarly to the fibre-like case the Euclidean coordinates of the oriented unit tangent vector $\mathbf{t}_i^j$ of the oriented geodesic curves 
$g(E_0, A_i^j)$ $((i,j)=(2,3),(3,2),$ $(1,2),(1,3))$ is given by (4.6).
Finally, we get similarly to the fibre-like case using the methods of spherical $\mathbf{S}^2$ geometry the following
\begin{lem}
The sum of the interior angles of any hyperbolic-like right angled hyperbolic-like right angled geodesic triangle is less or equal to $\pi$. 
\end{lem}
In the following table we summarize some numerical data of geodesic triangles for given parameters:    
\medbreak
\centerline{\vbox{
\halign{\strut\vrule\quad \hfil $#$ \hfil\quad\vrule
&\quad \hfil $#$ \hfil\quad\vrule &\quad \hfil $#$ \hfil\quad\vrule &\quad \hfil $#$ \hfil\quad\vrule &\quad \hfil $#$ \hfil\quad\vrule &\quad \hfil $#$ \hfil\quad\vrule
\cr
\noalign{\hrule}
\multispan6{\strut\vrule\hfill\bf Table 4, ~ ~ $y_2=1/2$ \hfill\vrule}%
\cr
\noalign{\hrule}
\noalign{\vskip2pt}
\noalign{\hrule}
z^3 & |\alpha_2^3|=|\alpha_3^2| & d(A_2A_3) & \omega_1 & \omega_2  & \sum_{i=1}^3(\omega_i)  \cr
\noalign{\hrule}
\rightarrow 0 & \rightarrow 0  & \begin{gathered} \rightarrow \\ \mathrm{arctanh}(1/2) \approx \\
\approx 0.5493 \end{gathered} & \rightarrow 0 & \rightarrow \pi/2 & \rightarrow \pi \cr
\noalign{\hrule}
1/10 & 0.0811 & 0.5638 & 0.1334 & 1.2830 & 2.9872 \cr
\noalign{\hrule}
1/3 & 0.2103 & 0.6994 & 0.3613 & 0.7170 & 2.6491 \cr
\noalign{\hrule}
\frac{999}{1000} & 0.0649 & 4.0707 & 0.5817 & 0.0913 & 2.2438\cr
\noalign{\hrule}
\frac{(10^6-1)}{10^6} & 0.0330 & 7.5174 & 0.0467 & 0.6112 & 2.2288\cr
\noalign{\hrule}
}}}
\medbreak
\subsubsection{Geodesic triangles with interior angle sum $\pi$}
In the above sections we discussed the fibre- and hyperbolic-like geodesic triangles and proved that there are right angled
geodesic triangles whose angle sum $\sum_{i=1}^{3}(\omega_i)$ is greater or equal to $\pi$, less or equal to $\pi$, but $\sum_{i=1}^{3}(\omega_i)=\pi$ is
realized if one of the vertices of a geodesic triangle $A_1A_2A_3$ tends to the infinity (see Table 3-4). We prove the following 
\begin{lem}
There is geodesic triangle $A_1A_2A_3$ with interior angle sum $\pi$ where its vertices are {\it proper} (i.e. $A_i \in \mathrm{SL}_2(\mathbb{R})$).   
\end{lem}
{\bf{Proof:}} We consider a hyperbolic-like geodesic right angled triangle with vertices $A_1=E_0=(1:0:0:0)$, $A_2=(1:0:y^2:0)$, 
$A_3^h=(1:0:0:z^3)$ and a fibre-like right angled geodesic triangle with vertices $A_1=E_0=(1:0:0:0)$, $A_2=(1:0:y^2:0)$, $A_3^f=(1:x^3:0:0)$
($0 < y^2,z^3 <1$, $0 < x^3 < \infty$). We consider the straight line segment (in Euclidean sense) $A_3^f A_3^h \subset \mathrm{SL}_2(\mathbb{R})$.

We consider a geodesic right angled triangle $A_1A_2A_3(t)$ where $A_3 \in A_3^f A_3^h$, $(t\in [0,1])$. $A_3$ is moving on the 
segment $A_3^hA_3^f$ and if $t=0$ then $A_3(0)=A_3^h$, if $t=1$ then $A_3(1)=A_3^f$. 

Similarly to the above cases the interior angles of the geodesic triangle $A_1A_2$ $A_3(t)$ are denoted by $\omega_i(t)$ $(i \in \{1,2,3\})$. 
The angle sum $\sum_{i=1}^{3}(\omega_i(0)) < \pi$ and $\sum_{i=1}^{3}(\omega_i(1)) > \pi$. Moreover the angles $\omega_i(t)$ change 
continuously if the parameter $t$ run in the interval $[0,1]$. Therefore there is a $t_E\in (0,1)$ where $\sum_{i=1}^{3}(\omega_i(t_E)) = \pi$. 
~ ~ $\square$

We obtain by the Lemmas of this Section the following
\begin{thm}
The sum of the interior angles of a geodesic triangle of $\SLR$ space can be greater, less or equal to $\pi$. 
\end{thm}
\subsection{Translation triangles}
We consider $3$ points $A_1$, $A_2$, $A_3$ in the projective model of $\SLR$ space (see Section 2). 
The {\it part of the translation curve} $a_k$ between the points $A_i$ and $A_j$  
$(i<j,~i,j,k \in \{1,2,3\}, k \ne i,j$ are called sides of the {\it translation triangle} $T_t$ with vertices $A_1$, $A_2$, $A_3$. 
We have seen in the Section 2 by the equations of the translation curves (see Table 2) that the translation curves are straight lines 
in the projective model. It is easy to see that the images of the translation curves by translations $\bT^{-1}$ are also straight lines 
because the translation is a collineation (see \ref{2.6}).  

Considering a translation triangle $A_1A_2A_3$ we can assume by the homogeneity of the $\SLR$ geometry that one of its vertex 
coincide with the origin $A_1=E_0=(1:0:0:0)$ and the other two vertices are $A_2(1:x^2:y^2:z^2)$ and $A_3(1:x^3:y^3:z^3)$. 

We will consider the {\it interior angles} of translation triangles that are denoted at the vertex $A_i$ by $\omega_i$ $(i\in\{1,2,3\})$.
We note here that the angle of two intersecting translation curves depend on the orientation of their tangent vectors. 

Similarly to the geodesic cases, in order to determine the interior angle sum $\sum_{i=1}^3(\omega_i)$ 
of translation triangle $A_1A_2A_3$ 
we define a \emph{translations} $\bT_{A_i}$, $(i\in\{2,3\})$ (see \ref{2.6}) that 
maps the origin $E_0$ onto $A_i$.  

We have seen that $\omega_1$ and the angle of geodesic curves with common point at the origin $E_0$ is the same as the 
Euclidean one therefore it can be determined by usual Euclidean sense.   

The parametrization of a translation curve in the model is given by the geographical sphere coordinates 
$(\lambda, \alpha)$, as longitude and altitude, $(-\pi < \lambda \le \pi, ~ -\frac{\pi}{2}\le \alpha \le \frac{\pi}{2})$, 
and the arc-length parameter $s \geq 0$ (see Table 2). 

We denote the oriented unit tangent vector of the oriented translation curves $(E_0, A_i^j)$ with $\mathbf{t}_i^j$ where
$(i,j)=(2,3),(3,2),(1,2),(1,3)$.
The Euclidean coordinates of $\mathbf{t}_i^j$ are the same as the Euclidean coordinates of points $A_i^j$. 

\begin{figure}[!ht]
\centering
\includegraphics[width=0.5\textwidth]{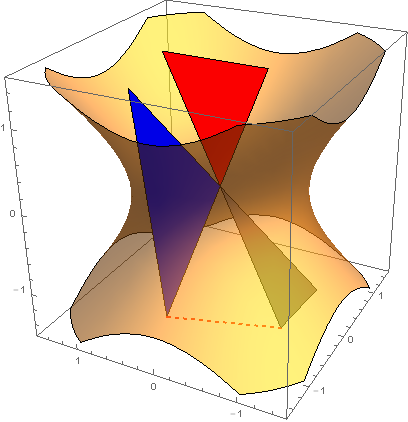}\includegraphics[width=0.5\textwidth]{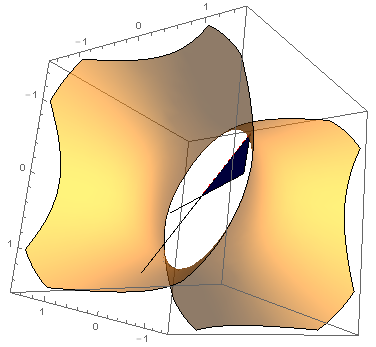}
\caption{Translations of the original (red) translation triangle.}
\label{fig:harom3}
\end{figure}

It was observed that the neighborhood of the origin  behaves like the Euclidean space so that the angle of two 
oriented tangent vectors with the origin as base point seems real size in the model and it can be determined by 
usual formula of the Euclidean geometry (see \ref{3.1}).
Now we translate both vertices $A_i$ by translations $\bT^{-1}_{A_i}$ to the origin to determine the other two interior angle 
of the translation triangle $A_1A_2A_3$:
\begin{equation}
\begin{gathered}
A_2=(1:x^2:y^2:z^2), ~ \bT^{-1}_{A_2}(A_2)=A_2^2=(1:0:0:0), {\bT^{-1}_{A_2}(A_1)}=A_1^2= \\ {=(1:-x^2:-y^2:-z^2)},
~ {\bT^{-1}_{A_2}(A_3)} = A_3^2=\Big( 1: 
 {\frac{x^2-x^3-y^2 z^3+y^3 z^2}{-x^2 x^3+y^2 y^3+z^2z^3-1}:}\\
 {\frac{-x^2 z^3+x^3 z^2+y^2-y^3}{-x^2 x^3+y^2 y^3+z^2z^3-1}:}
 { \frac{x^2 y^3-x^3 y^2+z^2-z^3}{-x^2 x^3+y^2 y^3+z^2z^3-1}\Big)}
\end{gathered}
\end{equation}
\begin{equation}
\begin{gathered}
A_3 =(1:x^3:y^3:z^3), ~ \bT^{-1}_{A_3}(A_3)=A_3^3=(1:0:0:0), ~ \bT_{A_3}^{-1}(A_1) = A_1^3 = \\ 
=(1:-x^3:-y^3:-z^3),~ {\bT^{-1}_{A_3}(A_2)} = A_2^3= \Big(1:
\frac{x^3-x^2-y^3 z^2+y^2 z^3}{-x^2 x^3+y^2 y^3+z^2z^3-1}: \\
\frac{-x^3 z^2+x^2 z^3+y^3-y^2}{-x^2 x^3+y^2 y^3+z^2z^3-1}:
\frac{x^3 y^2-x^2 y^3+z^3-z^2}{-x^2 x^3+y^2 y^3+z^2z^3-1} \Big).
\end{gathered}
\end{equation}
It is easy to see that the point pairs $(A_2,\bT^{-1}_{A_2}(A_1))$, $(A_3,\bT^{-1}_{A_3}(A_1))$ and $(\bT^{-1}_{A_2}(A_3),$ 
$\bT^{-1}_{A_3}(A_2))$ are antipodal.  
In Figure \ref{fig:harom3} there can be seen the translated triangles in the hyperboloid model. 
Also during the translation, the plane containing the triangle {\it twists, i.e.} the translated plane does not coincides generally 
with the original plane.
\begin{lem}
Let $\sigma=A_1A_2A_3$ be a plane in Euclidean sense trough the origin $(A_1=E_0)$ and be its Euclidean normal vector $\bv$. 
Then $\sigma$ is invariant for $\SLR$ translation $E_0Q$ where $Q \in \sigma$ if and only if $\bv$ is light--like (see ...). 
\end{lem}
\textbf{Proof:} 
The Euclidean equation of the plane $\sigma$ is:
\begin{equation}
\sigma:\ (y^2z^3-y^3z^2)x+(x^3z^2-x^2z^3)y+(x^2y^3-x^3y^2)z=0
\label{sikszig}
\end{equation}
The inhomogeneous coordinates of $\bT^{-1}_{A_2}(A_3)$ satisfies the (\ref{sikszig}) equation if and only if it is on $\sigma$:
\begin{equation}
\begin{gathered}
\frac{1}{-x_1 x_2+y_1 y_2+z_1z_2-1}\cdot
\left( (x_1-x_2)(y_1z_2-y_2z_1)-(y_1z_2-y_2z_1)^2+\right. \\  
\left.(y_1-y_2)(x_2z_1-x_1z_2)+
+(x_2z_1-x_1z_2)^2+(z_1-z_2)(x_1y_2-x_2y_1)+ \right.\\ \left. (x_1y_2-x_2y_1)^2\right)=
-(y_1z_2-y_2z_1)^2+(x_2z_1-x_1z_2)^2+(x_1y_2-x_2y_1)^2=\\ = -\bv_1^2+\bv_2^2+\bv_3^2=0
\Leftrightarrow \bv \mathrm{\ is\ light-like.} ~ ~ \square
\end{gathered} \notag
\end{equation}
Now we can claim the following theorem:
\begin{thm}
The sum of the interior angles of the translation triangle is greater or equal to $\pi$.
\end{thm} 
\begin{figure}[ht]
\centering
\includegraphics[width=6.5cm]{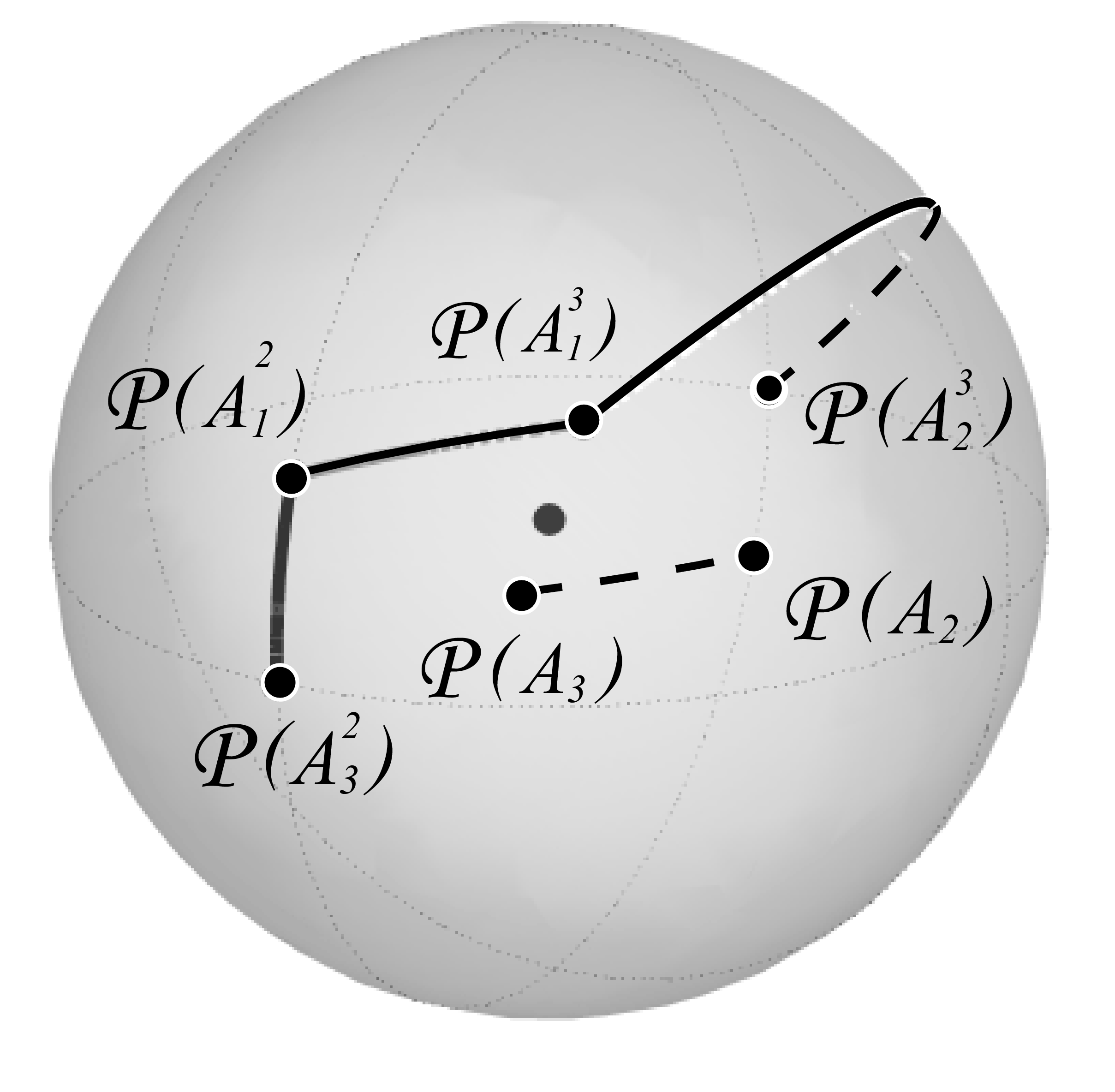}
\caption{Antipodal point pairs of the translation triangle projected onto the unit sphere}
\label{fig:antipod1}
\end{figure}
\textbf{Proof:} The translations $\bT_{A_2}^{-1}$ and $\bT_{A_3}^{-1}$ are isometries 
in $\SLR$ geometry thus $\omega_2$ is equal to the angle $((A_2^2 A_1^2), (A_2^2 A_3^2)) \angle$ (see Figure \ref{fig:antipod1}) 
of the oriented translation segments (Euclidean segments as well) $(A_2^2 A_1^2)$, $(A_2^2A_3^2)$ and $\omega_3$ is equal to the angle 
$((A_3^3 A_1^3),(A_3^3 A_2^3)) \angle$ 
of the oriented translation segments (Euclidean segments as well) $(A_3^3 A_1^3)$ and $(A_3^3 A_2^3)$ $(E_0=A_2^2=A_3^3$ see \ref{4.4}, \ref{4.5}).  

To get the angles $\omega_i$ we apply the $\cP$ projection from the origin onto the unit sphere around $E_0$ 
to the vertexes $A_2$, $A_3$ and $A_i^j$, $(i,j)=(1,2), (1,3), (2,3), (3,2)$. 
The measure of angle $\omega_i$ $(i\in \{1,2,3\})$ is equal to the spherical distance of the corresponding projected points on the unit sphere 
(see Figure \ref{fig:antipod1}). 
Due to the antipodality $\omega_1=\cP(A_2)A_1\cP(A_3) \angle =\cP(A_1^2)A_1 \cP(A_1^3) \angle$, therefore their corresponding spherical 
distances are equal, as well (see Figure \ref{fig:antipod1}). 
Now, the sum of the interior angles can be considered as three consecutive arc $(\cP(A_3^2) \cP(A_1^2))$, $(\cP(A_1^2) \cP(A_1^3))$, 
$(\cP(A_1^3) \cP(A_2^3))$ and $\cP(A_3^2)$
is antipodal to $\cP(A_2^3)$. 

Since the triangle inequality holds on the sphere, the sum of these arc lengths is greater or equal to the half 
of the circumference of the main circle on the unit sphere \textit{i.e.} $\pi$. 

By the Lemma 4.5 we obtain that the three consecutive arcs $(\cP(A_3^2) \cP(A_1^2))$, $(\cP(A_1^2) \cP(A_1^3))$, 
$(\cP(A_1^3) \cP(A_2^3))$ lie on a  main circle of the unit sphere if and only if the normal vector of the 
plane containing the triangle is light--like thus in this case the sum of the interior angles of the translation triangle is $\pi$. ~ ~ $\square$
%
%\begin{figure}[!ht]
%\centering
%\includegraphics[width=0.6\textwidth]{antipod1.png}
%\caption{Antipodal point pairs with main circle arcs.}
%\label{fig:antipod1}
%\end{figure}
%
%{\bf{Acknowledgement:}}

\end{document}